RESEARCH ARTICLE

# Transit facility allocation: Hybrid quantum-classical optimization

Einar Gabbassov  *

Department of Mathematics, Simon Fraser University, Burnaby, BC, Canada

* egabbassov@uwaterloo.ca

## Abstract

An essential consideration in urban transit facility planning is service efficiency and accessibility. Previous research has shown that reducing the number of facilities along a route may increase efficiency but decrease accessibility. Striking a balance between these two is a critical consideration in transit planning. Transit facility consolidation is a cost-effective way to improve the quality of service by strategically determining the desirable allocation of a limited number of facilities. This paper develops an optimization framework that integrates Geographical Information systems (GIS), decision-making analysis, and quantum technologies for addressing the problem of facility consolidation. Our proposed framework includes a novel mathematical model that captures non-linear interactions between facilities and surrounding demand nodes, inter-facility competition, ridership demand and spatial coverage. The developed model can harness the power of quantum effects such as superposition and quantum tunnelling and enables transportation planners to utilize the most recent hardware solutions such as quantum and digital annealers, coherent Ising Machines and gate-based universal quantum computers. This study presents a real-world application of the framework to the public transit facility redundancy problem in the British Columbia Vancouver metropolitan area. We demonstrate the effectiveness of our framework by reducing the number of facilities by 40% while maintaining the same service accessibility. Additionally, we showcase the ability of the proposed mathematical model to take advantage of quantum annealing and classical optimization techniques.

**Data Availability Statement:** All the data and the source code used to produce the experiments are uploaded to the GitHub repository: https://github.com/lostdevfound/bus-stop-optimization. Since the data is tightly integrated with the source code the GitHub repository is the best option.

**Funding:** The author received no specific funding for this work.

**Competing interests:** The authors have declared that no competing interests exist.

## 1 Introduction

In a growing society, public transportation systems are widely used to serve large populations in a space-time-concentrated travel demand manner [1]. The more successful a public transportation system is, the greater the positive effect on the population who uses it [2]. The success of a transportation system requires a balance between the accessibility and efficiency of the system [3]. Historically, some transit systems have issues achieving this balance.

Accessibility refers to how easily demand nodes (e.g. urban areas) can access the transit facilities [4]. Generally, the closer a facility is to a demand node, the more accessible the facility becomes. Thus, more facilities in place along a given route would create more accessibility for transit riders [4]. Efficiency is defined by how much distance a transit route can cover in a particular time window [4]. The efficiency of a transit route depends on a few factors such as the number of facilities on the route, traffic flow, dwell times at facilities, and the number of passengers boarding and alighting at each transit facility. Ideally, a transit route is highly efficient to appease customers and costs in travel time. However, finding a balance between accessibility and efficiency is preferred. As a case study, we consider the transportation system of one of the densest cities in British Columbia, Vancouver. Translink is a mass transportation system responsible for the regional transportation network of Metro Vancouver, including public transport, major roads and bridges. We use publicly available GIS data from Translink and Census program to optimize one of the most problematic bus routes. We demonstrate that the proposed framework can reduce the number of transit facilities by 40% while maintaining the same accessibility and improving efficiency.

### 1.1 Contributions

With these definitions for accessibility and efficiency, the main goal is to determine the best allocation of a limited number of transit facilities. This paper proposes a novel end-to-end optimization framework that accomplishes the goal. Specifically, the paper's contributions are:

- We propose a novel mathematical model inspired by the non-linear Spatial Interaction Coverage (SIC) models [5, 6]. The proposed model addresses transit facility allocation and accessibility by modelling population ridership (demand), physical accessibility (distances between demand nodes and facilities), connectivity (number of destinations from a facility), distance decay and inter-facility competition. One of the model's unique features is that its solutions could be naturally integral even if the integrality constraints are disregarded. See Section 2 and Section 4.

- Although expressive, SIC models are extremely hard to solve [5, 6]. We address this issue by using Quadratic Unconstrained Binary Optimisation (QUBO), which is closely related to the Ising spin glass model used in Statistical Physics and Quantum Simulation [7, 8]. Due to this relation, our model is readily solvable on a wide variety of classical and quantum optimization metaheuristics such as Quantum Annealing, quantum-inspired Digital Annealing [9], quantum optics-based hardware such as Coherent Ising Machines [10] and universal gate quantum computers [11]. Thus, the proposed model enables transportation planners to utilize and experiment with quantum-based technology along with a conventional classical-computing toolset. See Section 3.

- To analyze the quality of attained solutions, we present a theoretical upper bound on optimal values of the proposed model. The upper bound is fairly tight and can be computed in polynomial time. See Section 5.

- We propose a flexible integration of GIS into the model using the multi-criteria decision-making analysis [12]. The analysis is widely used in many applications related to decision-making in logistics, industry and business. See Section 3.3.

In addition to the above, we present a hybrid quantum-classical solver built with the open-source *dwave-hybrid* Python framework [13]. To explore solution search space robustly, the solver utilizes classical hill-climbing mechanics together with quantum superposition and quantum tunnelling effects [14].

The remainder of this paper proceeds as follows. Section 1.2 reviews various optimization models and algorithms. Section 2 describes the formulation of the proposed mathematical model. Section 3 describes QUBO formalism, model simplification, application of TOPSIS and quantum-classical solver. Sections 4 and Section 5 discuss the model's properties, such as solutions' integrality and an objective function bound. Section 6 showcases the application of the framework to the transit bus route in Vancouver, British Columbia. Section 8 summarizes both mathematical and numerical findings and highlights the implications of this study. Section 9 contains concluding remarks and future research directions.

## 1.2 Related literature

There are several modelling approaches for locating optimal facilities. However, these models do not simultaneously address distance decay, coverage range, competitors and GIS data. A well-known approach for addressing a service coverage range is Location Set Covering Problem (LSCP) [15]. LSCP seeks to place the minimum number of facilities such that all demand nodes can be served. While LSCP can identify the minimal subset of facilities necessary for covering all demand nodes it is unable to account for distance attenuation [5]. In some application scenarios, the model's requirement for all demands to be covered can be a fairly stringent condition that limits the choice of possible facility configurations. Moreover, often we would like to control the number of chosen facilities and not necessarily aim for the smallest number of them. Another modelling method is a distance-based approach such as the *p*-median problem [16]. The *p*-median formulation aims to minimize the total demand-weighted distance between each demand node and the nearest facility. While the *p*-median problem implements weighted distance as one of the decision factors and allows to control the number of chosen facilities, the binary nature of facility assignment does not allow for a demand node to be served by an arbitrary number of facilities.

Furthermore, recent studies actively integrate GIS into transit optimization [17]. GIS introduces more refined and data-oriented modelling, allowing us to account for census geographies such as population density per area, dissemination area boundaries and global positioning data. One of the significant requirements of GIS approach is large-scale data handling. Studies that incorporate GIS often use closed-source and inaccessible software which imposes strict requirements on how data should be handled. This in turn makes the proposed methodology not applicable to a wider range of scenarios. One such GIS-based approach is introduced in [3, 5]. The study combines the aspects of the *p*-median problem with the Maximal Covering Location Problem (MCLP) discussed in [5]. The proposed Spatial Interaction Coverage (SIC) model introduces the concept of interaction between a demand node and a facility. The interaction incorporates GIS-based demand, distance decay, geographical coverage and inter-facility competition. The model aims to maximize the interaction between a demand node and a facility. The resulting formulation is a non-linear, binary integer problem. The downside of the approach is the necessity to develop a proprietary optimization heuristic that needs tight coupling with commercial GIS software. The non-linearity of the problem combined with integrality constraints significantly reduces the choices of available solvers.

Several studies investigate an integrated Multiple Criteria Decision Making (MCDM) approach to optimize the facility-location problem. Such decisions require an algorithm that optimizes the quantitative and qualitative factors of each facility. These take the form of a scoring/ranking algorithm that takes into consideration factors including: a transit facility's connectedness to adjoining routes, ridership demand and its location [18]. For a comprehensive multi-criteria decision analysis [18] proposes Analytical Hierarchical Process (AHP) as such a method. AHP is a structured technique for organizing and analyzing complex decisions based on pairwise comparisons of desired criteria. AHP has several drawbacks: judgment criteria are not guaranteed to be always consistent and it does not scale well to large datasets [19]. An alternative to AHP is the Technique for Order Preference by Similarity to Ideal Solution (TOPSIS) [12, 20]. TOPSIS is based on an aggregating function representing "closeness to the ideal", which originated in the compromise programming method. The idea behind TOPSIS is computing a synthetic ideal candidate and determining the closeness of each facility to the ideal. Unlike AHP, TOPSIS scales well to large datasets, it can process hundreds of candidates and dozens of decision criteria, finding an optimal compromise in difficult decision situations.

Several well-known options exist for implementing and solving combinatorial problems. Commercial solvers, such as Lingo, Gurobi or CPLEX can be used, but also programming languages like Python allow the development of custom solvers that can tackle arbitrary combinatorial problems. Using objected-oriented programming has made it easier to integrate GIS data into a model formulation [3]. As mentioned previously, our paper discusses the formulation of a mathematical model using QUBO as a framework. In that regard, there exists a wide variety of powerful software- and hardware-based QUBO solvers which are successfully applied to large-scale optimization problems. QUBO problems can be efficiently solved on classical central processing units (CPU), graphics processing units (GPU) as well as specialized optimization hardware such as D-Wave's Quantum Annealer (DA) [21, 22], Fujitsu's Digital Annealer (DA) [9] and the Coherent Ising Machine (CIM) [10]. In addition to hardware solutions, there exist a number of software-based optimization algorithms that can efficiently solve QUBO problems. One of such algorithms is Tabu Search (TS). TS is a metaheuristic search method. Unlike local search methods which often get trapped in suboptimal regions, TS avoids this behaviour by prohibiting already visited solutions [23]. Simulated annealing (SA) is another probabilistic hill-climbing heuristic that searches for a global optimum [24]. Combining software and specialized hardware solutions into a single algorithm yields a powerful and versatile approach for solving optimization problems. We summarize both quantum and classical heuristics that can readily solve our proposed model in Table 1.

## 2 Model formulation

We base our formulation on the aforementioned MCLP and SIC models. To reflect the current demand along the route, we extensively utilize GIS Census Program data and work on the scale of the Dissemination Area (also known as the Census Block or Dissemination Block), which is one of the minimal geographic units in the population Census Program in North America. The objective of the facility allocation optimization problem is to select $p$ optimal facilities such that these $p$ facilities maximize route accessibility while reducing redundancy.

**Table 1. The summary of quantum and classical algorithms that support the proposed mathematical model.**

| Algorithms table | | |
| --- | --- | --- |
| **Name** | **Type** | **Operational principle** |
| QA [13, 21] | Quantum hardware | Quantum adiabatic evolution |
| CIM [10] | Quantum-optical hardware | Quantum-to-classical transition |
| QAOA [25] | Quantum gate algorithm | Quantum adiabatic evolution |
| QMF [26] | Quantum gate algorithm | Amplitude amplification by quantum interference |
| VQE [27] | Quantum gate algorithm | Variational Quantum State Preparation |
| DA [9] | Classical hardware | Markov Chain Monte Carlo search |
| SA [13, 24] | Classical algorithm | Monte Carlo, hill-climbing search |
| TS [13, 23] | Classical algorithm | Monte Carlo, hill-climbing search |
| PI-GNN [28] | Classical algorithm | Solution preparation with Graph Neural Networks |
| QANA [29] | Classical algorithm | Differential evolution search |
| APOPT [30, 31] | Classical algorithm | Sequential Quadratic Programming |
| IPOPT [31, 32] | Classical algorithm | Interior Point Method |

https://doi.org/10.1371/journal.pone.0274632.t001

We formulate the model as a maximization of a QUBO problem. All QUBO problems consist of *linear* and *quadratic* parts. As we will see later the *linear* part of our model accounts for relations between facilities and neighbouring demand nodes (Dissemination Areas) within a radius $R_0$. And the *quadratic* part accounts for inter-facility competition, i.e. if a facility has neighbouring competitors within a defined radius $R_1$ the facility has less importance during an optimization process. This approach helps to determine optimally spaced facilities that cover the maximum possible demand along the route. Our mathematical notation is summarized in Table 2.

We commence by introducing the model formulation in a quadratic form

$$\text{maximize} \sum_{j \in F} Q_{jj} x_j + \sum_{j \in F} \sum_{i \in F_j} Q_{ij} x_i x_j \quad (1)$$

subject to the constraint

$$\sum_{j \in F} x_j = p, \quad (2)$$

where

$$Q_{jj} = \frac{1 + 2m_j}{(1 + m_j)^2} w_j^\alpha \sum_{k \in A_j} a_k d_{jk}^{-\beta}, \quad (3)$$

$$Q_{ji} = \begin{cases} -\dfrac{Q_{jj}}{(1 + 2m_j)}, & \text{if } i \in F_j, \\ 0, & \text{otherwise.} \end{cases} \quad (4)$$

Our goal is to maximize the quadratic function (1) by activating an optimal combination of the $x_j$ for $j = 1, 2, \ldots, |F|$ such that there are a total of $p$ active facilities, i.e., the $p$-constraint (2) is satisfied. In other words, we want to select the best possible combination of $p$ facilities that maximize the objective. To develop an intuition behind (1) we look at its first *linear* part $\Sigma_{j \in F} Q_{jj} x_j$. Let

$$\tilde{Q}_{jj} = w_j^\alpha \sum_{k \in A_j} a_k d_{jk}^{-\beta}. \quad (5)$$

**Table 2. Notation for all variables used in the model.**

| | |
|---|---|
| $j$ | index of a facility |
| $F$ | set of all facilities along a route |
| $F_j$ | set of competitor facilities in a $R_1$-neighbourhood of a facility $j$ |
| $m_j$ | maximum number of possible competitors in $F_j$ |
| $A_j$ | set of demand nodes in a $R_0$-neighbourhood of a facility $j$ |
| $a_k$ | population at a demand node $k$ |
| $w_j$ | weight of a facility $j$ |
| $\alpha$ | exponent controlling weight |
| $d_{kj}$ | distance between a facility $j$ and a demand node $k$ |
| $\beta$ | exponent controlling distance $d_{kj}$ |
| $p$ | number of facilities to retain |
| $x_j$ | binary decision variable for a facility $j$. $x_j = 1$ if a facility is active |

https://doi.org/10.1371/journal.pone.0274632.t002

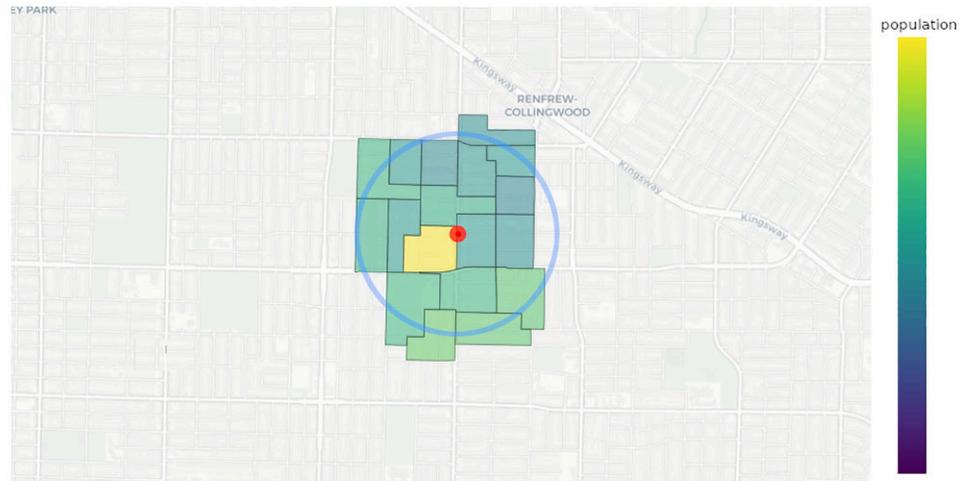

**Fig 1. An illustration of an urban map where the facility *j* (red dot) is surrounded by its neighbouring demand nodes (coloured urban blocks).** The blue circle determines the extent of the facility's neighbourhood $A_j$. The demand nodes are coloured according to their population density.

https://doi.org/10.1371/journal.pone.0274632.g001

This quantity represents the weighted sum of population $a_k$ attenuated by distance $d_{jk}$ between the facility *j* and the demand node *k*. We note that a demand node *k* is in the neighbourhood $A_j$ of a facility *j*. Thus each facility *j* has an associated aggregated population, see Fig 1. The greater the weighted population the higher the likelihood that the facility will be selected as a part of the final solution. We now explain the meaning of $\frac{1+2m_j}{(1+m_j)^2}$ in (3). As mentioned previously, our model accounts for inter-facility competition. We set the following simple rule: if a facility *j* has one neighbouring competitor then the aggregated population can be distributed among two facilities. Thus we divide $\tilde{Q}_{jj}$ by 2. Similarly, if a facility *j* has 2 competitors, the aggregated population can be distributed between 3 facilities, hence we divide $\tilde{Q}_{jj}$ by 3 and so on. This relation can be represented as follows:

$$f_j = \frac{\tilde{Q}_{jj}}{1 + \sum_{i \in F_j} x_i} x_j, \quad (6)$$

where $F_j$ is a neighbourhood of competing facilities. For example, if all the $x_i$ are active, the aggregated population of a facility *j* is divided by $1 + \sum_{i \in F_j} x_i = 1 + m_j$. However, the optimization problem based on (6) i.e.

$$\begin{aligned} \text{maximize} & \sum_{j \in F} f_j \\ \text{subject to} & \sum_{j \in F} x_j = p \end{aligned} \quad (7)$$

is non-linear and hard to solve. Therefore, we approximate it via a quadratic form by linearizing around the point $m_j$. To this end, put $z = \sum_{i \in F_j} x_i$. Then

$$\begin{aligned} f_j(z) &= \frac{\tilde{Q}_{jj}}{1+z} x_j \\ &\approx f_j(m_j) + f'_j(m_j)(z - m_j) \\ &= \frac{1 + 2m_j}{(1 + m_j)^2} \tilde{Q}_{jj} x_j - \frac{\tilde{Q}_{jj} x_j z}{(1 + m_j)^2} \\ &= \frac{1 + 2m_j}{(1 + m_j)^2} \tilde{Q}_{jj} x_j - \frac{\tilde{Q}_{jj}}{(1 + m_j)^2} \sum_{i \in F_j} x_i x_j. \end{aligned} \quad (8)$$

Finally, we observe that this gives the coefficients of the objective function

$$\frac{1 + 2m_j}{(1 + m_j)^2} \tilde{Q}_{jj} = Q_{jj}, \quad (9)$$

$$-\frac{\tilde{Q}_{jj}}{(1 + m_j)^2} = -\frac{Q_{jj}}{(1 + 2m_j)} = Q_{ji}. \quad (10)$$

This shows that all factors involving $m_j$ come from linearization around $m_j$ of (6) and that the *quadratic* part of the objective function models competition among the facilities.

Such a formulation has a number of interesting properties:

- The optimization problem can be intuitively studied from a graph-theoretical perspective (Section 3.2).

- The formulation can be easily cast into QUBO formalism.

- Both discrete and continuous optimization methods yield solutions to QUBO that are naturally binary, i.e. under certain conditions we can disregard integrality requirements and still obtain a binary solution (Section 4).

- We can derive a fairly tight upper bound for a maximum value of the objective function (Section 5).

### 2.1 Small example

We now look at a small scale example with 3 facilities, i.e. $F = \{1, 2, 3\}$. Without loss of generality, consider the second facility, with $j = 2$. Suppose that the facilities $j = 1$ and $j = 3$ are competitors of the second facility, thus the neighbourhood $F_2 = \{1, 3\}$. Suppose that $F_1 = \{2\}$ and $F_3 = \{2\}$. Noting that $Q_{ij} < 0$ for $i \neq j$ we have

$$\tilde{f}_1 = Q_{11} x_1 + Q_{21} x_2 x_1, \quad (11)$$

$$\tilde{f}_2 = Q_{22} x_2 + Q_{21} x_2 x_1 + Q_{23} x_2 x_3, \quad (12)$$

$$\tilde{f}_3 = Q_{33} x_3 + Q_{23} x_2 x_3. \quad (13)$$

And the objective function is

$$\tilde{f}_1 + \tilde{f}_2 + \tilde{f}_3. \quad (14)$$

For instance, deactivating the second facility ($x_2 = 0$) completely removes contributions of a facility $j = 2$ to the objective function. That is, $\tilde{f}_2 = 0$, while both $\tilde{f}_1$ and $\tilde{f}_3$ have a greater contribution to the objective function value because $Q_{21}x_2x_1$, $Q_{23}x_2x_3 = 0$. Whereas, deactivating $x_1$ will increase the importance of the facility $j = 2$. It is important to note that the competition relation is reflexive, i.e. if the facility $j = 2$ has the competitor $j = 3$ then $F_3$ also contains the facility $j = 2$. Such reflexive relation leads to a combinatorial explosion of possible configurations of $x_j$.

## 3 Optimization methods

This section describes QUBO formalism and how to use it to obtain the Ising equivalent model. Furthermore, we demonstrate simplification of the $p$-constraint that reduces the complexity of the problem and provides more controlled optimization results. We conclude the section with the GIS integration using TOPSIS and the presentation of the Hybrid quantum-classical solver.

### 3.1 QUBO formalism

So far we have presented a quadratic objective function with one explicit $p$-constraint (2). In order to obtain the final QUBO formulation, we need to incorporate the constraint into the objective. Specifically, for a positive scalar $\gamma$, we subtract a quadratic penalty

$$\gamma \left( \sum_{j \in F} x_j - p \right)^2 \quad (15)$$

from the objective function. Hence, we now consider the problem

$$\text{maximize} \sum_{j \in F} Q_{jj} x_j + \sum_{j \in F} \sum_{i \in F_j} Q_{ij} x_i x_j - \gamma \left( \sum_{j \in F} x_j - p \right)^2. \quad (16)$$

We note that the $p$-constraint (2) is now a soft constraint, i.e. it is incorporated into the objective function as a quadratic penalty. It can be shown that (16) can be compactly expressed in a matrix form. Noting that $x_i^2 = x_i$ for $i = 1, \ldots, |F|$, the quadratic objective can be rewritten in a matrix form as

$$\sum_{j \in F} Q_{jj} x_j + \sum_{j \in F} \sum_{i \in F_j} Q_{ij} x_i x_j = x^t Q x, \quad (17)$$

where $x \in \{0, 1\}^{|F|}$ and the matrix $Q \in \mathbb{R}^{|F| \times |F|}$ has entries $Q_{ij}$ defined in (3) and (4). Similarly, the penalty term can be expressed as

$$\left( \sum_{j \in F} x_j - p \right)^2 = (\mathbb{1}^t x - p)^t (\mathbb{1}^t x - p) = x^t \mathbb{1} \mathbb{1}^t x - 2p \mathbb{1}^t x + p^t p. \quad (18)$$

We can drop the constant term $p^t p = p^2$, and noting that $x_i^2 = x_i$ we can rewrite the above as

$$x^t \mathbb{1} \mathbb{1}^t x - 2p \mathbb{1}^t x = x^t \mathbb{1} \mathbb{1}^t x - 2p x^t x = x^t (\mathbb{1} \mathbb{1}^t - 2pI) x. \quad (19)$$

Let $G = \mathbb{1} \mathbb{1}^t - 2pI$. Then the final QUBO formulation is

$$x^t Q x - \gamma x^t G x = x^t (Q - \gamma G) x = x^t \bar{Q} x, \quad (20)$$

where $\bar{Q} := Q - \gamma G$ is a symmetric matrix. Therefore, our objective is to solve the following optimization problem

$$\max_{x \in \{0,1\}^{|F|}} x^t \bar{Q} x. \tag{21}$$

## 3.2 A graph theoretical view and partitioning the constraint

One of the advantages of the QUBO formulation is that we can view the entire problem as a network graph and assess its complexity in a visual manner. To this end, we view each non-zero $Q_{ij}$ of the quadratic objective presented in (1) as an edge of a graph with $|F|$ vertices. Since (1) relates pairs of facilities by their neighbourhood-dependent competition, we expect the resulting graph to be sparse: namely, each vertex (facility) is connected to a neighbouring vertex on a transit route. See Fig 2, left. Such a graph is almost a tree and it is easy to work with. However, upon the addition of the $p$-constraint (2) the resulting graph is a complete graph—see Fig 2, center. This can be intuitively explained by the fact that in order to satisfy the constraint each vertex must be aware of the activity status ($x_j = 1$ or $x_j = 0$) of every other vertex. This implies that the decision variable such as $x_1$ is related to $|F| - 1$ other decision variables. To mitigate the complexity, we split the $p$-constraint into several constraints

$$\sum_{j \in S_k} x_j = p_k \text{ for } k = 1, \ldots, m \tag{22}$$

where $\{S_k\}_{k=1}^m$ is a partition of the facility set $F$ and all the $p_k$ add up to $p$. Then, the quadratic penalty in (15) is replaced by

$$\gamma \sum_{k=1}^m \left( \sum_{j \in S_k} x_j - p_k \right)^2. \tag{23}$$

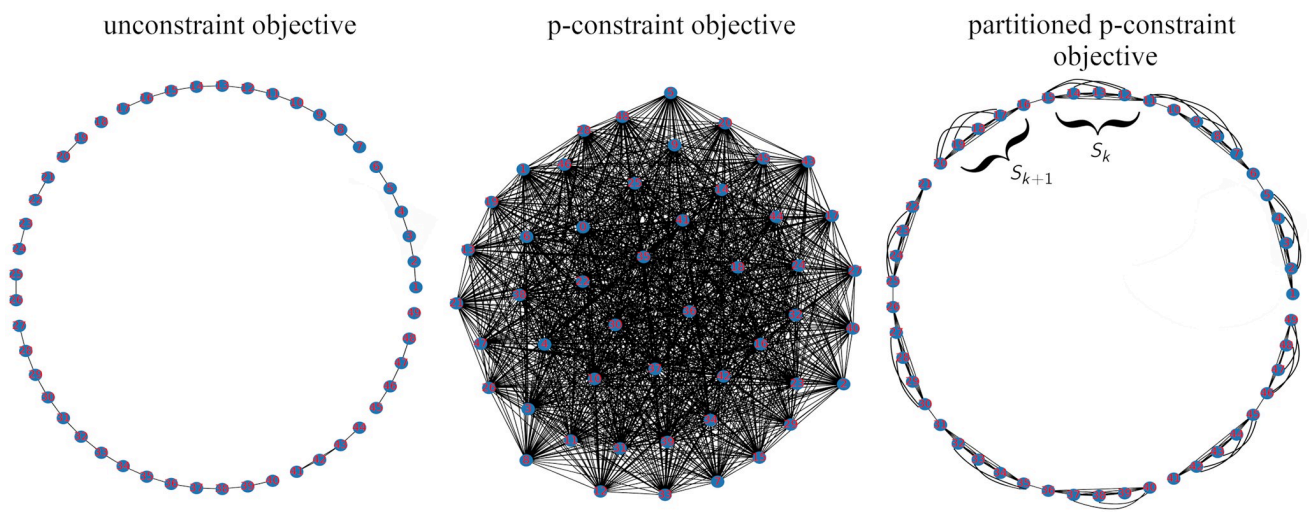

**Fig 2. Graph networks depicting relations between 49 facilities under three different objective functions.** Each vertex $j \in \{1, \ldots, 49\}$ represents the facility $j$, and each edge $(i, j)$ represent competition between the facilities $i$ and $j$. (Left) Inter-facility relations under the unconstrained quadratic objective (1). (Center) Objective function with the $p$-constraint added as a penalty (16). The resulting graph is dense because each facility must be aware of all other facilities. (Right) The objective function with the partitioned $p$-constraint as a penalty (24). The resulting graph is sparse because each facility must be aware of only facilities in the subset $S_k$ to which it belongs and its neighbouring competitors.

https://doi.org/10.1371/journal.pone.0274632.g002

The new objective function is

$$\sum_{j \in F} Q_{jj} x_j + \sum_{j \in F} \sum_{i \in F_j} Q_{ij} x_i x_j - \gamma \sum_{k=1}^{m} \left( \sum_{j \in S_k} x_j - p_k \right)^2. \quad (24)$$

For instance, we can partition a bus route in batches of 5 facilities and require every batch to have 4 active facilities at all times. This gives $S_1$ = {1, 2, 3, 4, 5} with $p_1$ = 4, $S_2$ = {6,7,8,9,10} with $p_2$ = 4 and so on. For the last subset $S_m$ with at most 5 remaining facilities we set $p_m$ = min{$|S_m|$, 4}. Just like the original constraint, the partitioned constraint enforces $p$ active facilities at all times. However, it has one important advantage: the resulting graph is sparse, see Fig 2 right. In this instance, the degree of $x_1$ is at most 4. Another important advantage of this approach is that it ensures that feasible solutions do not have multiple neighbouring facilities removed from a route. Such undesirable removal could result in substantial gaps along a route. Hence, this approach allows for more granular control over which facilities are removed. Finally, we can compactly express the new objective in (24) in a matrix form

$$\max_{x \in \{0,1\}^{|F|}} x^t \hat{Q} x, \quad (25)$$

where $\hat{Q}$ is a symmetric matrix.

### 3.3 Weights estimation, TOPSIS

The TOPSIS algorithm is widely used in multi-attribute decision-making problems such as supply chain logistics, engineering and facility allocation problems [12, 20]. As mentioned before, the algorithm determines the best ranking (weighting) for a set of candidates with certain attributes. The algorithm treats each candidate as a point in the Euclidean space, where the candidate's attributes become spatial coordinates. Each candidate is weighted based on its Euclidean distance to the synthesized ideal solution. The ideal solution is a theoretical candidate with the most desired attributes, and if such a candidate existed, it would get the highest possible ranking. In practice, the ideal candidate may not exist, but it could be synthesized based on the available data and criteria determined by a decision-maker. TOPSIS creates the ideal solution and produces the final ranking (weights) for all candidates based on their Euclidean distance proximity to the ideal solution.

In order to estimate the model's weights $w_j$ for $j$ = 1, .., $|F|$, we use TOPSIS from the open-source Python library *scikit-criteria* [33]. First, we create an *alternatives matrix* that consists of rows of attributes of each facility. The matrix serves as an input to the TOPSIS algorithm. In this study, we consider three attributes to determine the TOPSIS weights $w_j$. Namely, we use

1. Monthly demand data per facility (demand attribute).

2. Transit facility connectedness to other services (connectedness attribute).

3. Number of neighbouring landmarks that could produce additional demand (landmarks attribute). For example, the landmarks could be hospitals, schools, subways, malls etc.

We note that these data are independent of the parameters used in the optimization model. Table 3 is an example of such a matrix for a problem with three facilities. The $W_d$, $W_c$, $W_l$ are user-defined importance parameters of each attribute. In this case, we give the highest preference ($W_d$ = .45) to the demand attribute of each facility.

We summarize the TOPSIS workflow in the Fig 3. Note that despite the fact that the facility A has the highest demand, the Fig 3 suggests that facility C has the highest weighting as it has

**Table 3. A user constructed three-facility matrix for TOPSIS algorithm.** We prioritize a facility's demand attribute by setting $W_d = 0.45$. If the main goal was to retain a well-connected allocation then $W_c$ would have the largest weight.

| facility priority weights | demand | connectedness | landmarks |
|---|---|---|---|
| | $W_d = .45$ | $W_c = .3$ | $W_l = .25$ |
| facility A | 5064 | 1 | 7 |
| facility B | 4383 | 3 | 6 |
| facility C | 4735 | 3 | 8 |

https://doi.org/10.1371/journal.pone.0274632.t003

better overall characteristics. We compute weights $w_j$ by generalizing this method to all facilities.

### 3.4 Model optimization

Taking advantage of QUBO's equivalence to the transverse Ising model [22] we can view the optimization problem as energy minimization of the Ising objective function. Briefly, the transverse Ising model was originally developed to describe quantum phenomena in

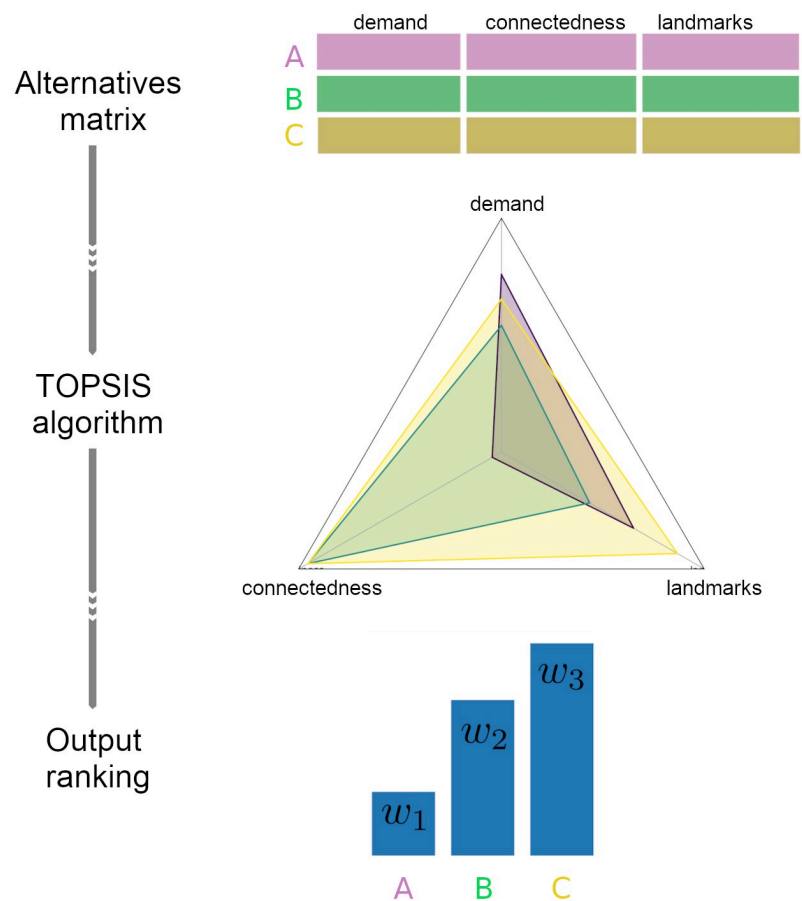

**Fig 3. A visual schematic of the weight estimation procedure with the TOPSIS algorithm.** First, the *alternatives matrix* is constructed and supplied to the TOPSIS algorithm. TOPSIS embeds the attributes into the Euclidean space. The three facility attributes are mapped on three axes: demand, connectedness and landmarks. From the radar plot, it is clear that facility C has overall better characteristics. Hence, it gets a higher ranking. The output is the facility weights.

https://doi.org/10.1371/journal.pone.0274632.g003

magnetism. The model consists of discrete variables that represent atomic spins that can be in one of two states (+ 1 or −1). It has been proposed that quantum annealing is an effective physical process that can attain the minimum energy state of a spin system [22, 34]. In 2007 a Canadian quantum computing company D-Wave Systems released a quantum device that utilized quantum annealing to solve Ising model-related problems. Since then they have released more powerful hardware suitable for large-scale discrete optimization problems. The D-Wave Quantum Processing Unit (QPU) can be viewed as a heuristic that efficiently minimizes Ising objective functions using a physically-realized version of quantum annealing. In this paper, we utilize D-Wave's quantum hardware along with classical computing heuristics to find the best solution for our problem.

### 3.5 Hybrid solver

It has been shown that decomposing QUBO problems into smaller QUBO sub-problems and combining the resulting solutions can be an effective method for obtaining optimal results [35, 36].

To find the best solution, we use an open-source *dwave-hybrid* Python framework [13]. The hybrid solver attempts to solve the entire problem at once while also solving multiple child sub-problems of the original problem. We integrate the classical heuristic methods tabu search and simulated annealing with quantum annealing into a single system. Both tabu search and simulated annealing work on the entire problem, whereas quantum annealing works on multiple child sub-problems of a smaller size where only high energy impact variables are optimized. The high energy impact variables are the variables whose contribution to the objective function is greater than the rest of the variables. All three heuristics run in parallel branches and combine their best solutions during each iteration $k$. See Fig 4. The process terminates whenever the time or iteration budget is depleted. We summarize the algorithm structure in Algorithm 1. The use of such a hybrid approach is highly beneficial for several reasons. First, we recall that tabu search and simulated annealing are hill-climbing local search heuristics. During the minimization process, local search heuristics tend to get trapped in suboptimal regions of an optimization landscape. To avoid being trapped in local minima, both methods deploy mechanisms of accepting uphill moves (worsening moves) that could lead to better

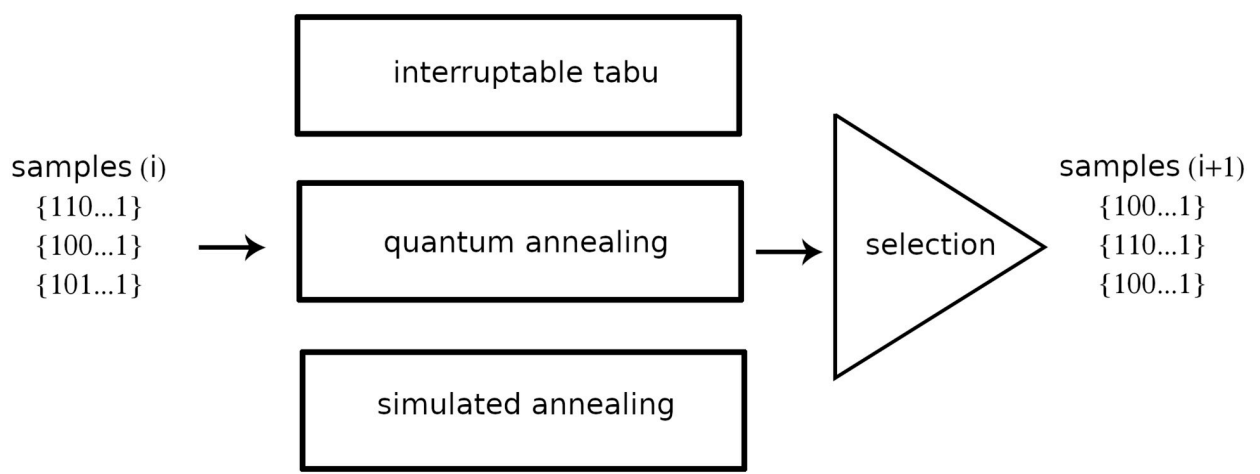

**Fig 4. A simplified schematic of a single iteration cycle of three branches of solvers.** The SA and TS solvers perform the classical hill-climbing optimization and provide redundancy in the case of decoherence errors in the quantum annealer. The quantum annealer extends the hill-climbing search with quantum tunnelling through optimization landscape barriers in the direction of the global minimum.

https://doi.org/10.1371/journal.pone.0274632.g004

regions. Tabu search can accept uphill moves based on its tabu list and fitness characteristics of a candidate move. By contrast, simulated annealing uses a probabilistic rule called the Metropolis-Hastings criterion [37] for accepting uphill moves. While hill-climbing turns out to be a fairly successful approach, it is not the only way of navigating the optimization landscape. Quantum annealing is also capable of hill-climbing, but it can also perform quantum tunnelling [14]. Intuitively, this can be seen as moving towards the global optimum right through an optimization landscape barrier. Such a tunnelling move is useful when the barrier is so tall that hill-climbing methods would need to accept too many suboptimal moves to traverse the barrier. Therefore, quantum tunnelling is a useful addition that aids the exploration of the solution search space.

**Algorithm 1** The hybrid algorithm runs simulated annealing (SA), tabu search (TS) and quantum annealing (QA) on multiple sub-problems.

```
Require:
1: Problem matrix Q with dimensions n × n.
2: Maximum iterations K.
3: Maximum time T.
4: Sub-problem size m ≤ n. If m = n the sub-problem is the entire
problem.
5: Begin
6: Initialize a guess solution x*.
7: Set k ← 1.
8: while k ≤ K and run time <T do
9:    x_sa ← run SA on a problem x^tQx starting with x*.
10:   x_tb ← run TS on a problem x^tQx starting with x*.
11:   x* ← choose deterministically or stochastically x_sa or x_tb.
12:   Identify a subset X of size m such that for j ∈ X, x_j has high
energy impact.
13:   Formulate a sub-problem y^tQy such that y_i = x_i* is fixed for i ∉ X.
14:   x_qa ← run QA on a problem y^tQy where for j ∈ X, y_j is a variable.
15:   If needed, repeat Steps 12-14 for different X and track x_qa.
16:   x* ← choose the best of x_qa, x_sa, x_tb.
17:   k ← k + 1.
18: end while
19: Return x*.
20: End
```

One might wonder why not entirely rely on quantum computation which can do both hill climbing and tunnelling. The answer is fairly nuanced. In short, the current generation of quantum hardware suffers from noise and errors that occur during quantum computation [38]. Without error correction, the accumulation of errors can significantly degrade the quality of the optimization results and even distort the quantum implementation of a problem [39]. Therefore, relying exclusively on quantum computation is still not entirely feasible. The hybrid approach, on the other hand, takes the best of both worlds where the parallel execution branches create necessary redundancy. Should one system fail due to errors or get trapped, the other branch takes over and provides its solution for the next iteration.

## 4 Integrality of solutions under continuous relaxation

In this section, we investigate the properties of a matrix $\hat{Q}$ in (25) and demonstrate that under mild conditions, the optimization problems (21) and (25) naturally yield integral solutions even if the integrality constraints are disregarded. First, $\hat{Q}$ is a symmetric matrix. Due to the nature of the problem and our specific formulation choices, the matrix $\hat{Q}$ has relatively large positive diagonal entries. This can be observed by inspecting (2). Based on a numerical

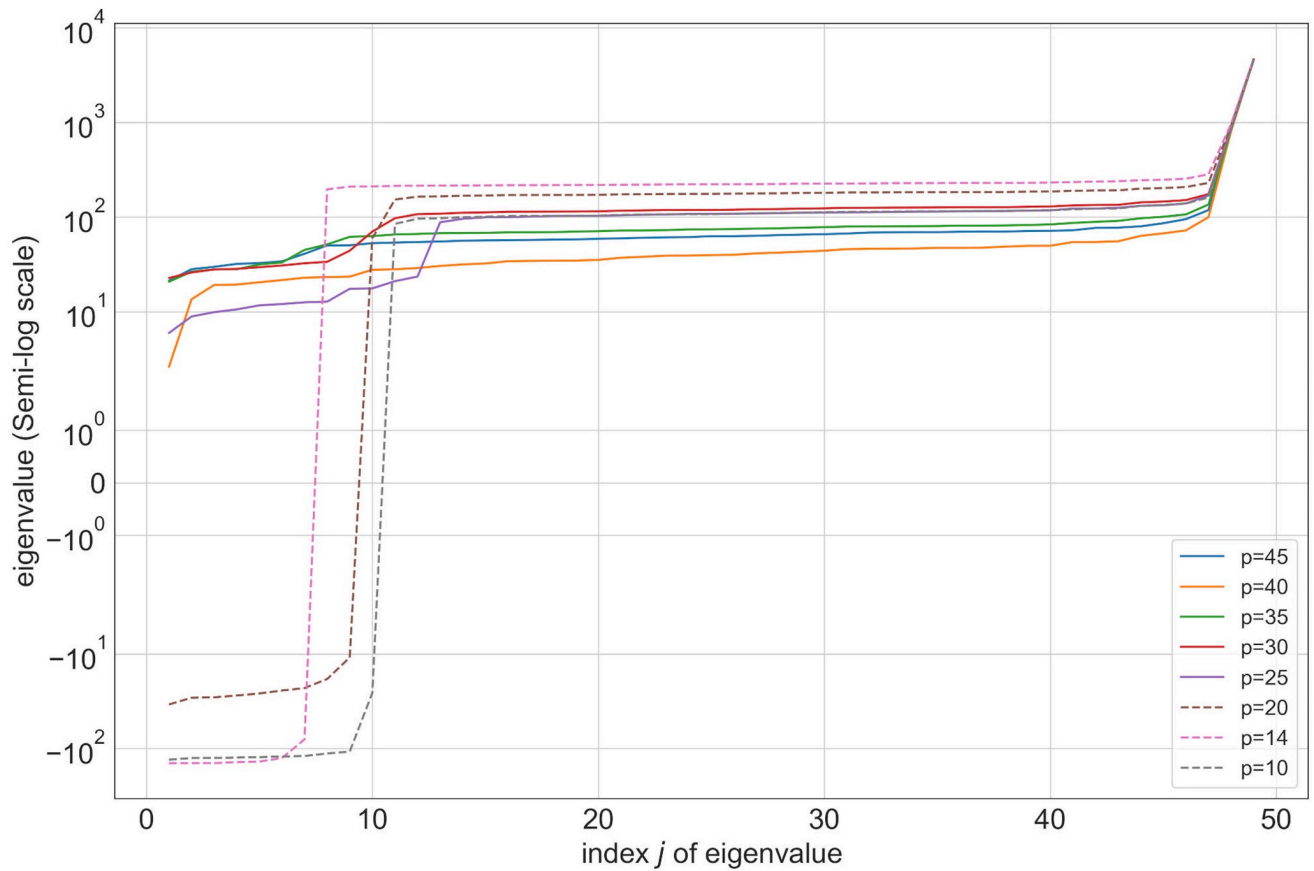

**Fig 5. Eigenvalues of $\hat{Q}$ for |F| = 49 and different $p$ values.** The objective function with $\frac{p}{|F|} \leq 0.5$ is nonconvex. The eigenvalues of $\hat{Q}$ for $p = 10$, $p = 14$, $p = 20$ are drawn with dotted lines to highlight the loss of positive definiteness and hence convexity.

https://doi.org/10.1371/journal.pone.0274632.g005

investigation of $\hat{Q}$ we learn that the matrix is positive definite whenever $\frac{p}{|F|} > 0.5$ and $m_j \leq 2$, see Fig 5. In other words, $\hat{Q}$ is positive definite whenever the ratio of retained facilities to a total number of facilities is more than a half and there are at most 2 competitors per facility. In the case when $\hat{Q}$ is not positive definite, we can obtain an equivalent representation of the problem with a perturbed positive definite matrix $Q' = \hat{Q}(\lambda + \epsilon)I$ where $\lambda$ is the smallest eigenvalue of $\hat{Q}$ and $\epsilon$ is chosen such that $\lambda + \epsilon > 0$ [40].

Since $\hat{Q}$ is positive definite, it is informative to consider the continuous relaxation of the problem:

$$\text{maximize } x^t \hat{Q} x$$
$$\text{subject to} \qquad\qquad\qquad\qquad (26)$$
$$0 \leq x_j \leq 1 \text{ for } j = 1, \ldots, |F|.$$

We relax the problem by allowing each $x_j$ to be a continuous variable restricted to the compact interval [0, 1]. The resulting constraint set is a hyper-cube. Hence, if $x$ is a feasible solution, then $x \in [0, 1]^{|F|}$. The positive definiteness of $\hat{Q}$ implies that the objective function in (26) is convex. Given that we aim to maximize the objective and $\hat{Q}$ is positive definite, the optimal solution occurs at one of the corners of the hyper-cube. Since the position of each corner corresponds to a binary string of length $|F|$, the optimal solution $x^*$ is guaranteed to be binary, i.e. $x_j^* \in \{0, 1\}$ for $j = 1, \ldots, |F|$. This means, if we allow $0 \leq x_j \leq 1$, then any suitable continuous optimization solver will yield binary solutions. It is worth noting that even though the objective function is convex we do maximization. Hence, there is no guarantee that the problem will be solved in polynomial time.

Using the continuous optimization solver IPOPT (Interior Point OPTimizer) along with the open-source GEKKO optimization framework [31] we obtain a solution to (26) that is identical to the solution to (25) from the Hybrid solver presented in the previous section. Furthermore, we do another solve using the branch and bound solver APOPT (Advanced Process OPTimizer) and obtain exactly the same solution again. Given the consensus between hybrid, continuous and combinatorial solvers we have a reason to believe that the obtained solution is globally optimal. Moreover, this demonstrates the flexibility of our formulation as it can be solved on a wide variety of discrete and continuous solvers with minimal modifications. In the next section, we investigate the optimality of the solution by presenting an upper bound for the problem.

## 5 Upper bound

We would like to ensure that the solutions to the problem are close to a theoretical upper bound for the problem. We claim that for all solutions $x \in \{0, 1\}^{|F|}$ that satisfy the partitioned $p$-constraint (22) the following holds,

$$x^t \hat{Q} x \leq \sum_{i=1}^{p} \lambda_i, \tag{27}$$

where $\lambda_1 \geq \ldots \geq \lambda_p \geq \ldots \geq \lambda_{|F|}$ are eigenvalues of $Q$ defined in (3) and (4) and $p$ is the number of facilities we want to retain. To show this, let $|F| = n$ and $\Sigma_k p_k = p$. Then

$$\begin{aligned} x^t \hat{Q} x &= \sum_{j \in F} Q_{jj} x_j + \sum_{j \in F} \sum_{i \in F_j} Q_{ij} x_i x_j - \gamma \sum_{k=1}^{m} \left( \sum_{j \in S_k} x_j - p_k \right)^2 \\ &= \sum_{j \in F} Q_{jj} x_j + \sum_{j \in F} \sum_{i \in F_j} Q_{ij} x_i x_j. \end{aligned} \tag{28}$$

We note that the quadratic penalty term is zero because $x$ satisfies the partitioned $p$-constraint, $\sum_{j \in S_k} x_j = p_k$ for all $k$. From the definition of $Q$ we know that for all $i$, $Q_{ii} > 0$ and $Q_{ij} \leq 0$ for $i \neq j$. Hence let us drop the off-diagonal terms,

$$\sum_{j \in F} Q_{jj} x_j + \sum_{j \in F} \sum_{i \in F_j} Q_{ij} x_i x_j \leq \sum_{j=1}^{n} Q_{jj} x_j. \tag{29}$$

Let $\{Q_{i_k i_k}\}_{k=1}^{n}$ be a non-increasing re-arrangement of $\{Q_{ii}\}_{i=1}^{n}$. Since the partitioned $p$-constraint (22) implies $p$-constraint (2), we have $\sum_{i=1}^{n} x_i = p$. This implies there are $n - p$ decision variables $x_i$ that are equal to zero. Hence,

$$\sum_{j=1}^{n} Q_{jj} x_j \leq \sum_{k=1}^{p} Q_{i_k i_k}. \tag{30}$$

Applying the Schur-Horn theorem [41, Theorem 4.3.45] yields

$$\sum_{k=1}^{p} Q_{i_k i_k} \leq \sum_{k=1}^{p} \lambda_k. \tag{31}$$

Therefore, for all $x \in \{0, 1\}^n$ satisfying (22), we have

$$x^t \hat{Q} x \leq \sum_{i=1}^{p} \lambda_i. \tag{32}$$

This shows that an upper bound can be calculated by summing the $p$ largest diagonal entries or eigenvalues of $Q$.

## 6 Case study

In this section, we apply our framework to the northbound bus route B20 of the transit system in Vancouver, British Columbia. The route connects suburban areas to Vancouver downtown and has 49 bus stops. The following sections present accessibility analysis and route efficiency improvements for transit-rider friendly value $p$ = 40 (20% reduction of bus stops). In addition, we perform a redundancy analysis and show that it is possible to remove up to 40% of the bust stops while maintaining the same accessibility.

### 6.1 Accessibility and efficiency analysis

Using our open-source GIS framework, we create a computer model of the northbound route and surrounding dissemination areas (demand nodes). The model stores information about the population and geographical vector boundaries of each dissemination area, bus stop connectivity to other routes, TOPSIS weights $w_j$, neighbourhoods $F_j$ and $A_j$ for each bus stop $j$. As mentioned before, the northbound route B20 has 49 bus stops. We set $p$ = 40, which is roughly 20% reduction of the number of stops. We use a radius of 400 meters for the neighbourhood $F_j$ (neighbourhood of competitors around the stop $j$) and a recommended 500 meters walking distance for the radius of $A_j$ (neighbourhood of dissemination areas around the stop $j$) with $\alpha$, $\beta$ = 1. Even though the number of bus stops is reduced by 20%, the retained stops still cover the same number of dissemination areas. This means the recommended walking distance to a bus stop of 500 meters is also retained. Using bus dwell time and boarding/alighting data from Translink we estimate a reduction of travel time by approximately 7 minutes, which is 13% of travel time on the unoptimized route during work hours. Therefore, we reduce the number of bus stops by 20% while preserving the same route accessibility and increasing service effectiveness. Fig 6 plots unoptimized (original) and optimized routes together with dissemination areas.

In Fig 7 we present the best objective function values found with the Hybrid solver for different values of $p$. Each best value is compared with its corresponding upper bound estimate presented in the previous section. The smallest ratio of the best solution to the upper bound is 5954/6291 ≈ 0.95 and it occurs at $p$ = 45. This implies that the best objective function values are fairly close to the unattainable theoretical upper bound. Moreover, for $p$ > 30, an increase in the number of bus stops does not significantly improve the effectiveness of the route. Hence we may conclude that any value close to $p$ = 35 will provide a good balance between efficiency and accessibility. For example, low values of $p$ will significantly impact accessibility. Indeed, for $p$ = 14 there is a visible loss in coverage, see Fig 8.

### 6.2 Redundancy analysis

To further analyze the best solutions for different $p$, we count how many northbound bus stops are within the walking distance of each dissemination area. The counting yields sequences of numbers $N_p = \{n_1^{(p)}, \ldots, n_i^{(p)}, \ldots, n_m^{(p)}\}$ for $p$ = 1, 2, ..., 49 where each entry $n_i^{(p)}$ corresponds to the number of bus stops within the walking distance of the dissemination area $i$

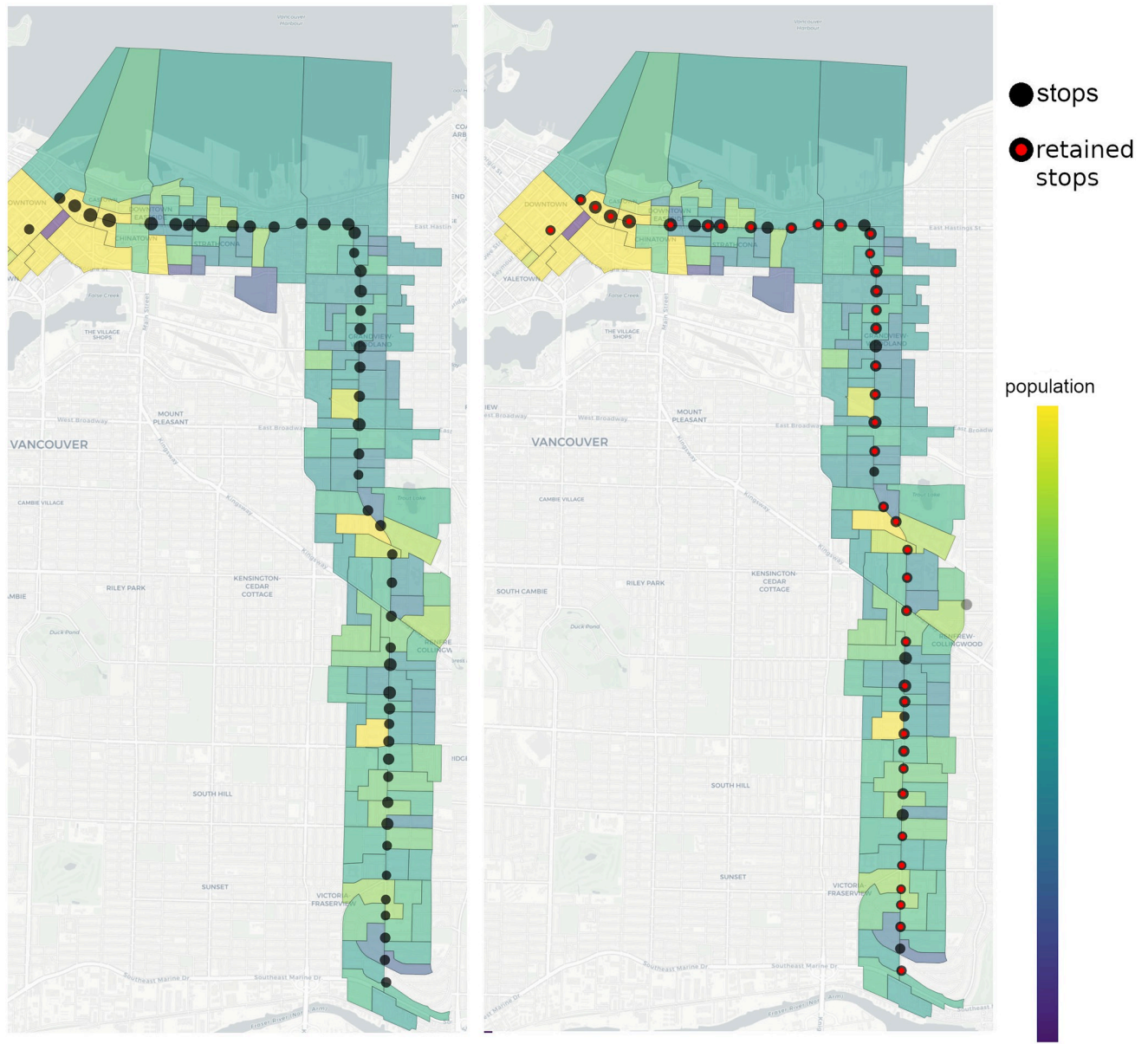

**Fig 6. Vancouver northbound B20 route with dissemination areas; unoptimized (Left), optimized (Right).** The black dots (bus stops) have different sizes according to their weight $w_j$. The plot on the right demonstrates the optimized route where the retained stops are marked in red.

https://doi.org/10.1371/journal.pone.0274632.g006

given that $p$ stops are retained. For example, for $p = 10$ we have $N_{10} = \{3, 0, \ldots, 1\}$. This means that the first, second and the last dissemination areas have 3, 0 and 1 bus stops within the walking distance respectively. For each $p$ we compute the minimum, maximum and median of the $N_p$. We summarize these statistics in Fig 9. Fig 9 demonstrates that for $p = 49$ every dissemination area has a minimum of 2 stops and some areas have up to 6 stops. This indicates redundancy. For $p = 30$ the best solution still covers every area (minimum is 1). Whereas at $p = 25$ we start to see the loss of coverage as there are areas with no stops within the walking distance. From this figure, we can conclude that even when 40% of stops are removed the model still manages to find an optimal configuration of retained stops to provide complete coverage. This analysis coincides with the marginal increase of the objective function for $p > 30$ (Fig 7, yellow line). Further, this suggests that the model effectively handles inter-facility competition.

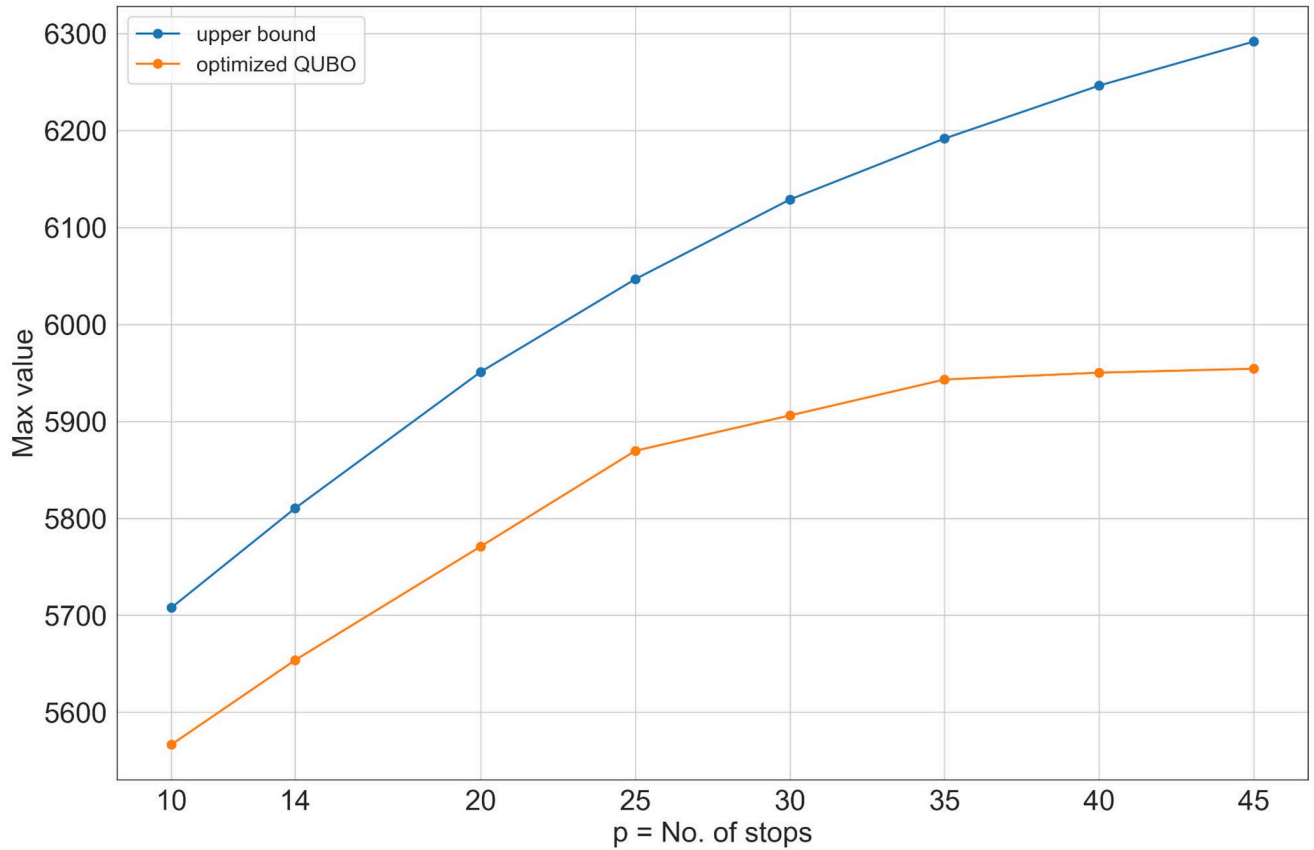

**Fig 7. The objective function maximum for *p* retained stops and corresponding theoretical upper bound.**

### 6.3 Additional numerical experiments

This section compares our model solved with the hybrid solver against a SIC type model solved with the APOPT solver. We demonstrate that the proposed methodology yields statistically better optimization results in terms of demand coverage for different numbers of retained facilities.

For this experiment, we randomly generate 30 synthetic problems with 50 facilities and 135 demand nodes. Each facility *j* has uniformly random weights $w_j \in (0, 1]$. The synthetic facility and demand node placements are based on the randomly perturbed coordinates of the facilities and demand nodes of the case study route B20. We choose the maximum perturbation offset of 150 meters because it allows generating diverse synthetic problems while still representing a relatively realistic scenario.

Each generated instance is solved using two approaches. The first approach $A_1$ uses the proposed mathematical formulation in (16) and the Hybrid solver presented in 3.5. We set one second time limit on the hybrid solver for each problem instance. The second approach $A_2$

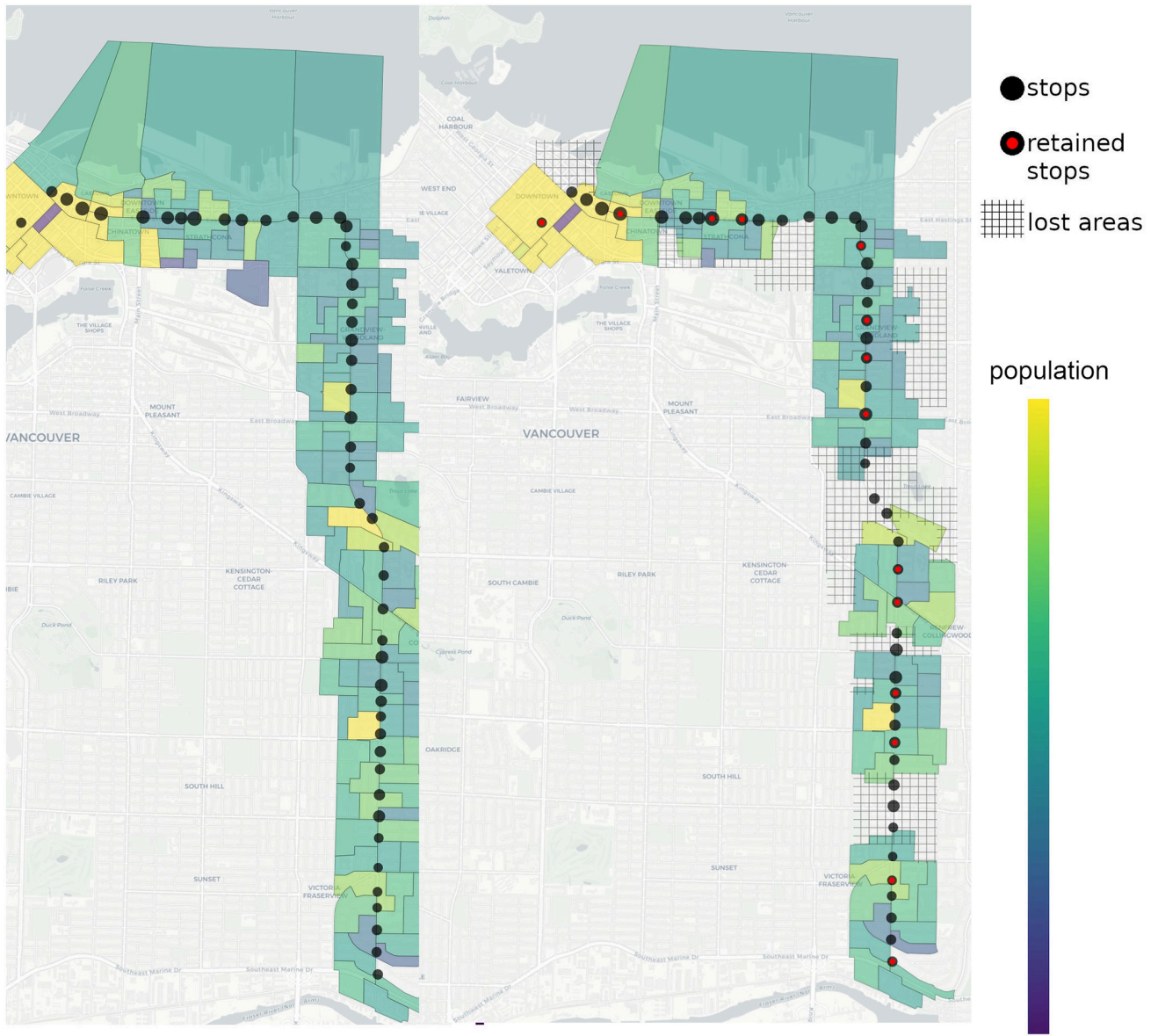

**Fig 8.** Unoptimized route (Left). Optimized route (Right). The best solution for *p* = 14 does not cover all the dissemination areas, resulting in an accessibility decrease.

https://doi.org/10.1371/journal.pone.0274632.g008

uses the SIC type formulation with the APOPT solver. The SIC formulation is as follows

$$\max \sum_{i \in D} \sum_{j \in F} S_{ij}$$

Subject to

$$S_{ij} = \frac{a_i w_j d_{ij}^{-1} x_j}{\sum_{k \in N_i} w_k d_{ik}^{-1} x_k + 1} \quad \text{for all } i, j \quad (33)$$

$$\sum_{j \in F} x_j = p$$

$$x_j \in \{0, 1\} \quad \text{for all } j.$$

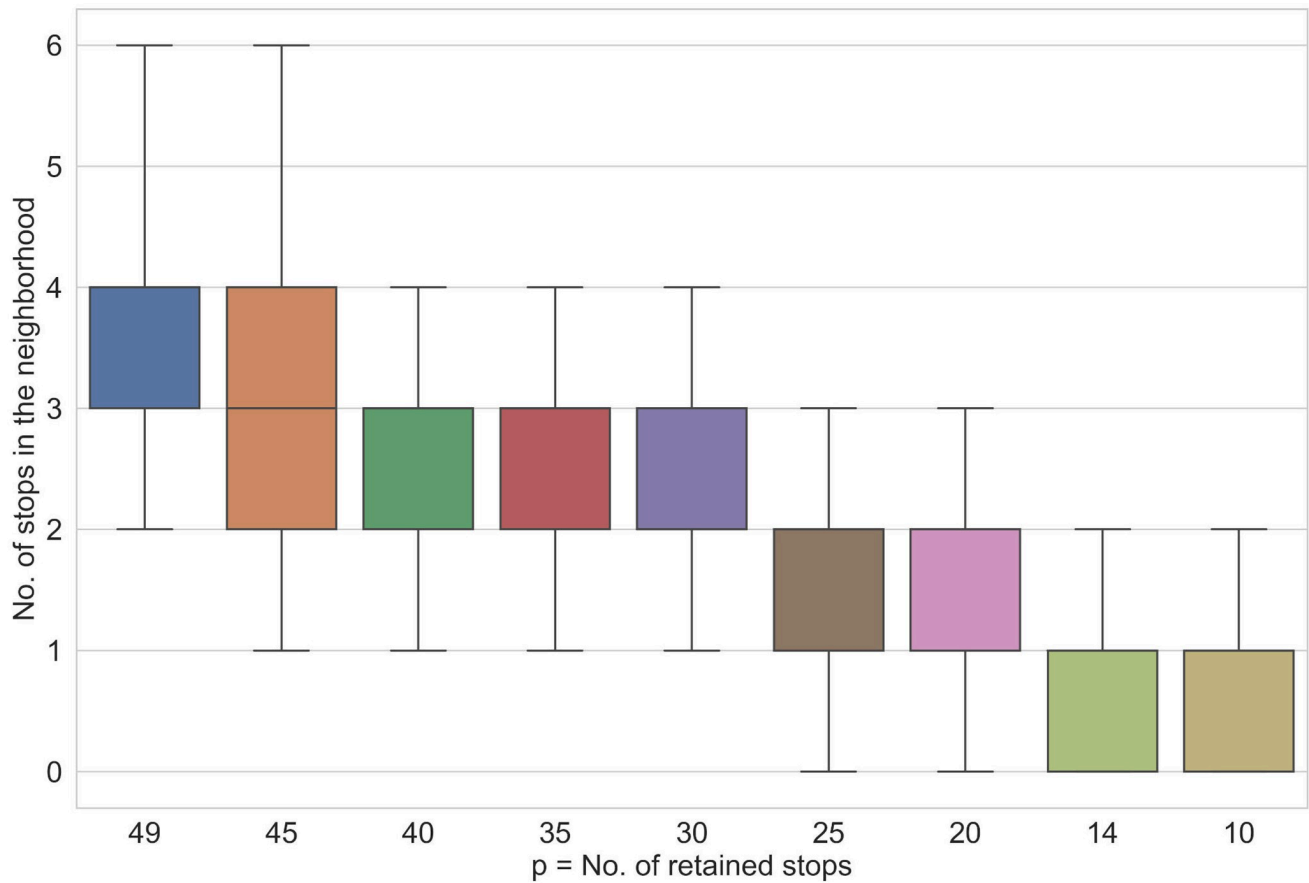

**Fig 9. Min, max and a median number of stops in dissemination area neighbourhoods for different values of *p*.**

https://doi.org/10.1371/journal.pone.0274632.g009

In this problem, $D$ denotes the set of all demand nodes, and $N_i$ denotes the set of facilities in the neighbourhood of the demand node $i$. The rest of the notation is defined as in Table 2. The SIC problem in (33) always has a solution for all values of $p \in \{1, \ldots, |F|\}$ and any location configuration of facilities and demand nodes. We set no time limit on the APOPT solver for each problem instance and allow 100000 maximum iterations to solve the problem. The convergence to a solution happens before the depletion of the iteration budget.

We test two hypotheses: the null hypothesis ($H_0$), which states the first approach $A_1$ covers less demand than $A_2$, and an alternative hypothesis ($H_a$) which states that $A_1$ covers more demand than $A_2$. The demand coverage difference $d$ is the difference between the number of demand nodes covered using $A_1$ and the number of demand nodes covered using $A_2$. If $d > 0$, then $A_1$ is better. Using the Wilcoxon signed-rank test [42] and the synthetic dataset, we summarize the results in Table 4.

## 7 Comparison with related models

This section qualitatively compares and summarizes the proposed model with the related models reviewed in Section 1.2. Specifically, we compare the proposed model with LSCP [15], $p$–median [16], MCLP and SIC models [5].

We commence the comparison by emphasizing the fact that the proposed model is primarily designed for quantum optimization. This means the proposed model is the only model

Table 4. *p*-values from the Wilcoxon signed-rank test on the 30 synthetic instances with the significance level $\alpha$ = 0.05.

| number of retained stops | *p*-value | reject $H_0$ | average *d* |
|---|---|---|---|
| 45 | 0.143 | No | 0.3 |
| 40 | 0.045 | Yes | 1.0 |
| 35 | 0.005 | Yes | 3.4 |
| 30 | 0.001 | Yes | 3.6 |
| 25 | 0.002 | Yes | 5.9 |

https://doi.org/10.1371/journal.pone.0274632.t004

amenable for quantum optimization and capable of taking advantage of quantum phenomena such as quantum superposition, quantum tunnelling and entanglement. Another unique differentiating feature of the proposed model is the ability to produce naturally integral solutions even if the integrality constraints are disregarded, and continuous optimization methods are applied. This is not the case for all other related models, as they require heuristic techniques that enforce the integrality of solutions. Finally, it is possible to compute a tight upper bound on the optimal objective function value. This computation can be done efficiently as it is only necessary to compute the sum of *p* largest diagonal entries of the matrix *Q*—for further details, see Section 5.

We summarize the comparison of the proposed model with the related models in Table 5.

## 8 Results and discussion

Unlike most literature on facility allocation and route optimization, we explore non-linear, physics-inspired mathematical modelling. This allows us to attain an expressive model that captures the most substantial modelling aspects of P-median, Maximal Covering Location and Spatial Interaction Coverage models. During the model derivation process, we demonstrate that it captures complex interactions between facilities and demand nodes, inter-facility competition and distance decay. We show that GIS and demand data can be flexibly integrated using the decision-making analysis techniques such as TOPSIS. To demonstrate our model's computational robustness, we implement a classical-quantum Hybrid solver that harnesses the quantum effects such as superposition and quantum tunnelling while providing redundancy with classical solvers. As a case study, we apply our optimization pipeline on the B20 bus route located in one of the most populated cities of Canada—Vancouver. We demonstrate that the proposed approach reduces the number of bus stop facilities by up to 40% while maintaining the same service coverage. Using the derived theoretical upper bound, we analyze the quality of solutions. The smallest ratio of the best solution to the upper bound is 5954/6291 ≈ 0.95,

Table 5. The summary of facility allocation models.

**Facility allocation models**

|  | LSCP | *p*–median | MCLP | SIC | Our model |
|---|---|---|---|---|---|
| Loss of coverage | No | No | Yes | No | Yes |
| Distance decay | No | No | No | Yes | Yes |
| Inter-facility competition | No | No | No | Yes | Yes |
| Coverage range | Yes | No | Yes | Yes | Yes |
| Naturally integral solutions | No | No | No | No | Yes |
| Quantum optimization | No | No | No | No | Yes |
| Polynomial time bounds | No | No | No | No | Yes |

https://doi.org/10.1371/journal.pone.0274632.t005

occurring at $p = 45$. This suggests that the best objective function values are fairly close to the unattainable theoretical upper bound. We conduct additional numerical experiments on the synthetically generated dataset and demonstrate that the proposed model solved with the Hybrid solver yields statistically superior demand coverage than a SIC class model solved with the APOPT solver.

The coming end of Moore's law [43] suggests that traditional computer devices are struggling to support the growing computational requirements. Due to this, research in alternative computational methods has steadily increased. Quantum computation holds a promise to speed up some of the most demanding optimization processes. This study clearly shows that the proposed optimization framework allows facility planners to consider emerging alternatives to classical computation by experimenting with novel specialized hardware that is capable of harnessing quantum effects.

While the ability to exploit novel quantum phenomena is one of the motivating factors for developing the proposed framework, the most crucial criterion is retaining important modelling aspects present in the predecessor models such as MCLP and SIC. One of the improvements present in our model is the possibility of introducing convexity into the objective function. We argue that whenever the objective function is convex, solutions to the problem are always binary, even if the integrality constraints are disregarded. This allows the proposed model to be treated from the perspective of continuous optimization that is rich in theory and optimization methods.

GIS and demand data can be large and complex [44]. Extracting useful information and deciding what aspects of the data should contribute to the decision-making process might be a great challenge [44, 45]. We address this by introducing decision-making analysis which allows facility planners to rank important factors in a formal analytical manner. Using TOPSIS, multiple competing decision criteria can be readily integrated into the optimization model for further analysis.

Although successful in many aspects, our approach suffers from the following limitations:

- As the number of retained facilities decreases, there is an inevitable loss in coverage for some demand nodes. When the demand coverage is critically important, and the solutions that lose coverage are deemed infeasible, our model fails as it allows for loss of coverage.

- Due to QUBO formalism, the addition of inequality constraints results in extra auxiliary binary variables. Possibly some inequality constraints may impact an optimization landscape's ruggedness, challenging the hill-climbing heuristics.

- Introducing integer variables that are not binary is a challenge. While not impossible, the result is likely to be inefficient.

We now discuss these limitations. While the loss of coverage may be considered a limitation, the ability to lose coverage guarantees the existence of a solution for all values of $p \in \{1, 2, \ldots, |F|\}$ and any geographical configuration of facilities and demand nodes. It may be impossible to reduce the number of facilities in many scenarios without losing some demand coverage.

Additional inequality constraints require auxiliary binary variables and need to be incorporated into an objective function. While this may pose a challenge, the non-linear nature of the model allows expressing many complex relations in an alternative fashion; for various methods see [46, 47]. Nevertheless, complex inequality constraints may pose mathematical and optimization issues.

The use of integer variables is another challenge from the quantum hardware perspective. The current generation of quantum hardware does not have enough resources to solve non-binary problems. However, recent developments [48] address this issue through various software-based methods.

## 9 Conclusion and future works

This paper presents a new framework for optimizing facility allocation problems. The framework aims to reduce the number of facilities while maximizing accessibility. It considers several factors such as geographical Census Program data, facility connectivity, demand data, inter-facility competition and proximity to public landmarks. We augment combinatorial optimization with multi-criterion decision-making analysis to establish data-driven rankings of each facility. This allows for building a parameterizable model that can be easily adapted for different usage scenarios.

The proposed mathematical model possesses unique characteristics. It has an efficient computable tight bound that can be used to analyze the quality of solutions. Under mild conditions, problem solutions are naturally integral. The model's equivalence to the Ising spin glass model allows harnessing quantum effects such as superposition and quantum tunnelling.

In future work, we plan to adapt the developed methods to multi-objective problems with distinct objectives capturing various optimization criteria. Since the proposed formalism allows algorithm implementation on fully programmable gate-based quantum computers, we plan to implement a hybrid optimization technique based on Variation Quantum Eigensolver [27] suitable for near-term quantum applications.

## Acknowledgments

We would like to thank professors Alexander Rutherford, Ben Adcock and Leonid Chindelevitch for their valuable help and advice. Also, we would like to thank Moe Antar for his help with GIS visualization and large dataset handling. Additional thanks to Sarah Morse for her initial writing contribution to the project. Many thanks to Parveen Sarana at Translink for sharing the ridership data. Special thanks to Dr. Mohammad Amin and Dr. Hossein Sadeghi at D-Wave for their advice on quantum computing. We are particularly grateful for the guidance given by Dr. Maliheh Aramon.

## Author Contributions

**Investigation:** Einar Gabbassov.

**Methodology:** Einar Gabbassov.

**Software:** Einar Gabbassov.

**Visualization:** Einar Gabbassov.

**Writing – original draft:** Einar Gabbassov.

**Writing – review & editing:** Einar Gabbassov.